\documentclass[CJK,12pt]{{article}}
\usepackage{}
\usepackage{mathrsfs}
\usepackage{amsfonts}
\usepackage{amsfonts,amssymb,amsmath,mathrsfs}
\usepackage{fancyhdr}
\usepackage{graphicx}
\usepackage{ifpdf}
\ifpdf
  \usepackage[CJKbookmarks]{hyperref}
\else
  \usepackage[CJKbookmarks,hypertex]{hyperref}
\usepackage{color,latexsym,amsfonts}
 \newcommand{\qed}{\hfill\rule{2mm}{3mm}\vspace{4mm}}

 \newtheorem{theorem}{Theorem}[section]
 \newtheorem{lemma}[theorem]{lemma}
 \newtheorem{corollary}[theorem]{Corollary}
 \newtheorem{proposition}[theorem]{Proposition}

 \def\beqlb{\begin{eqnarray}}\def\eeqlb{\end{eqnarray}}
 \def\beqnn{\begin{eqnarray*}}\def\eeqnn{\end{eqnarray*}}

 \def\proof{\noindent{\it Proof.~~}}
 \def\qed{\hfill$\Box$\medskip}
 \def\no{\nonumber}

\begin{document}

\bigskip

\noindent{\Large{\bf Branching structure for the transient random walk on a strip in a random environment \footnote{The project is partially supported by the National
Natural Science Foundation (Grant No. 11131003) of China.}}}

\bigskip\bigskip

 \centerline{Wenming Hong  \ and  \ Meijuan Zhang}

 \centerline{School of Mathematical Sciences, Beijing Normal
 University,}

 \centerline{Beijing 100875, People's Republic of China}

\smallskip

 \centerline{{\tt wmhong@bnu.edu.cn} and {\tt
 zhangmeijuan@mail.bnu.edu.cn}}

\bigskip\bigskip

\begin{center}
\begin{minipage}[c]{15cm}

\mbox{}\textbf{Abstract} An intrinsic branching structure within the transient random walk on a strip in  a random environment  is revealed. As  applications, which enables us to express the hitting time explicitly, and specifies the density of the absolutely continuous invariant measure for the ``environments viewed from the particle". \vspace{0.2cm}

\mbox{}\textbf{Keywords:}\quad Branching structure; random walk on a strip; random environment; hitting time, invariant measure, environments viewed from the particle.\\
\mbox{}\textbf{Mathematics Subject Classification}:  Primary 60J80;
secondary 60G50.
\end{minipage}
\end{center}


\section{Introduction}

\setcounter{equation}{0}

Let~$d\geq1$~be any integer and denote~$\mathscr{D}=\{1,2,\cdots,d\}$,~we consider random walks in a random environment on the strip~$S= \mathbb{Z} \times\{1,2,\cdots,d\}$. This model was introduced  by Bolthausen and Goldsheid (\cite{[BG00]}, 2000), where the conditions for recurrent and transient has been obtained. After that, Goldsheid (\cite{[Gol08]}, 2008) considered the hitting time of the walk by the method of ``enlarged random environments"; Bolthausen and Goldsheid (\cite{[BG08]}, 2008) obtained   the $(\log t)^2$ asymptotic behaviour and  Roitershtein  (\cite{[Roi08]}, 2008) proved a strong law of large numbers and an annealed central limit theorem for the walk in a suitable environment situation; etc..

The aim of this paper is to reveal the intrinsic branching structure within the transient random walk on a strip in  a random environment, which enables us to express the hitting time explicitly. Roitershtein (Theorem 2.3, \cite{[Roi08]}, 2008)  figured out the stationary distribution for the Markov chain of ``environments viewed from the particle" which is equivalent to the original distribution. To specify the density of the absolutely continuous invariant measure is another application of our branching structure. And as a by product, the rate of the LLN can be obtained.

For the nearest random walk in random environment (RWRE, for short) on the line, as we known,  the branching structure is a powerful tool in the proof of the famous result about ``stable law" by Kesten et al (\cite{[KKS75]}, 1975), and is also used by Ganterta and Shi  in (\cite{[GS02]}). The branching structure for the one dimensional RWRE with bounded jumps has been considered by Key (\cite{[Key87]}, 1987),  Hong \& Wang (\cite{[HW09]}, 2009) and Hong \& Zhang (\cite{[HZ10]}, 2010).

\subsection{Description of the model.}

We adapt the description of the model as that of \cite{[BG00]}.
Let $(P_n,Q_n,R_n), -\infty < n < \infty $, be a strictly stationary
ergodic sequence of triples of $m\times m$ matrices with non-negative elements such that for
all $n$ the sum $P_n + Q_n + R_n$ is a stochastic matrix, i.e., $(P_{n}+Q_{n}+R_{n})\mathbf{1}=\mathbf{1}$,~where~$\mathbf{1}$~is a column vector whose components are all equal to~$1$. We write the components
of $P_n$ as $P_n(i, j ), 1\leq i, j \leq m $, and similarly for $Q_n$ and $R_n$. Let ~$(\Omega,\mathscr{F},P,\theta)$
be the corresponding dynamical system with $\Omega$
 denoting the space of all sequences
$\omega := (\omega_n) = ((P_n,Q_n,R_n))$ of triples described above, $\mathscr{F}$ being the corresponding
natural $\sigma$-algebra, $P$ denoting the probability measure on $(\Omega,\mathscr{F})$, and the shift operator on~$\Omega$~defined by~$\theta$:~$(\theta \omega)_{n}=\omega_{n+1},~n\in \mathbb{Z}$.
The random walk on the strip~$S=\mathbb{Z} \times \mathscr{D}:= \mathbb{Z} \times\{1,2£¬\cdots, d\}$~is denoted by~$X=\{X_{n},n \in \mathbb{Z}\}$,
$$X_{n}=(\xi_{n},Y_{n}),~~\xi_{n}\in \mathbb{Z},~Y_{n}\in \mathscr{D}.$$
$\xi_{n}$~is the~$\mathbb{Z}$-coordinate of the walk and~$Y_{n}$~takes values in~$\mathscr{D}:=\{1,2£¬\cdots, d\}$.

For describing the initial distribution, we introduce  ~$\mathcal {M}_{d}$,
\begin{equation}
\mathcal {M}_{d}=\Big\{(\mu_{\omega})_{\omega\in \Omega}: \mu_{w}  \mbox{ is a probability measure vector on} \  \mathscr{D}=\{1,2,\cdots, d\}\Big\}.\no
\end{equation}

Given a environment~$\omega\in \Omega$~ and an  ~$\mu=(\mu_{\omega}) \in \mathcal {M}_{d}$, ~one can define the random walk~$X_{n}$~on the strip~$S= \mathbb{Z} \times \mathscr{D}$ to be a time-homogeneous Markov chain taking values in~$ \mathbb{Z} \times\{1,2,\cdots,d\}$, which is determined by its transition probabilities~$\mathfrak{Q}_{\omega}(z,z_{1})$:

\begin{displaymath}
\mathfrak{Q}(z,z_{1})
=\left\{\begin{array}{r@{\quad \quad}l}
P_{n}(i,j)  &  \mbox{if}~  \   z=(n,i),z_{1}=(n+1,j),
\\ R_{n}(i,j)  &  \mbox{if}~  \   z=(n,i),z_{1}=(n,j),
\\ Q_{n}(i,j)  &   \mbox{if}~  \  z=(n,i),z_{1}=(n-1,j),
\\ 0 & \mbox{otherwise},
 \end{array}\right.
\end{displaymath}
and initial distribution
\begin{equation}\label{id}
P_{\omega}^{\mu}(\xi_{0}=0,~Y_{0}=z_{0})=\mu_{\omega}(z_{0})\quad \mbox{for any}\quad z_{0}\in \mathscr{D}.
\end{equation}
This defines for any starting point~$x_{0}=(0,y_{0})\in S$~and for any $\omega \in (\Omega,\mathscr{F},P)$,~the quenched law~$P_{\omega}^{\mu}$~for the Markov chain by
\begin{equation}\label{tp}
P_{\omega}^{\mu}(X_{0}=x_{0},X_{1}=x_{1},\cdots,X_{n}=x_{n}):=\mu_{\omega}(y_{0}) \mathfrak{Q}_{\omega}(x,x_{1}) \mathfrak{Q}_{\omega}(x_{1},x_{2}) \cdots \mathfrak{Q}_{\omega}(x_{n-1},x_{n}).
\end{equation}
Then we define a annealed law~$\mathbb{P}^{\mu}=P\bigotimes P_{\omega}^{\mu}$~on~$(\Omega\times (\mathbb{Z} \times \mathscr{D})^{\mathbb{N}},\mathscr{F}\times \mathscr{G})$~by
\begin{equation}\label{ap}
\mathbb{P}^{\mu}(F\times G)=\int_{F}P_{\omega}^{\mu}(G)P(d\omega)\quad F\in \mathscr{F},~G \in \mathscr{G},
\end{equation}
and the expectation with respect to~$\mathbb{P}^{\mu}$~defined by~$\mathbb{E}^{\mu}$.~Statements involving~$P_{\omega}^{\mu}$~and~$\mathbb{P}^{\mu}$~are called quenched and annealed, respectively.

\


\noindent {\it Notations and assumption.}
Throughout the paper we use the notation~$\mathbf{0}=(0,0,\cdots,0)\in \mathbb{R}^{d}$,~$\mathbf{1}=(1,1,\cdots,1)\in \mathbb{R}^{d}$,~and denote ~$\mathbf{e_{i}}=(0,\cdots,1,\cdots,0),~(i=1,2,\cdots,d)$~as the canonical basis of~$\mathbb{R}^{d}$.
For the vector~$\mathbf{x}=(x_{j})$~and matrix~$A=(a(i,j))$,~define
$$\|\mathbf{x}\|:=\max_{j}|x_{j}| \quad\mbox{and}\quad \|A\|:=\max_{i}\sum_{j}|a(i,j)|.$$
We say that~$A$~is strictly positive~(denoted by~$A>0$)~if all its components satisfy~$a(i,j)>0$,~and~$A$~is non-negative~(which is denoted by~$A\geq0$)~if all~$a(i,j)$~are negative. If a~$d\times d$~real matrix~$A$~is non-negative,~$\|A\|:=\|A\mathbf{1}\|$.~Finally, we use the notation~$\mathbf{I}_{A}$~for the indicator function of the set~$A$.
For the random walk~$X_{n}=(\xi_{n},Y_{n})$,~we often use the expressions like~$\lim_{n\rightarrow\infty}X_{n}=+\infty$~which simply means~$\xi_{n}$~tends to~$+\infty$~as~$n\rightarrow\infty$.

The {\it hitting time} ~$T_{n}$~is defined as the the first time when the walk ~reaches layer~$n$~, $L_{n}:=\{(n,j),1\leq j\leq m\}$~  starting from a point~$z\in L_{0}:=\{(0,j),1\leq j\leq m\}$. Let $T_o=0$, and for  $n\geq 1$,
\begin{equation}\label{ht}
T_{n}:=\inf\{t:~X(t)\in L_{n}\} \quad \mbox{and} \quad \tau_{n}:=T_{n}-T_{n-1},
\end{equation}
with the usual convention that the infimum over an empty set is~$\infty$~and~$\infty-\infty=\infty$.\\

The following Condition C is  from   Bolthausen and Goldshied \cite{[BG00]}.

 \noindent{\textbf{Condition C.}}\\
C1 The dynamical system~$(\Omega,\mathscr{F},\mathbb{P},\mathcal {T})$~is ergodic.\\
C2\begin{equation}\label{c2}
\mathbb{E}\log(1-\| R_{n}+P_{n}\|)^{-1}<\infty \quad \mbox{and}\quad \mathbb{E}\log(1-\| R_{n}+Q_{n}\|)^{-1}<\infty.
\end{equation}
C3 For all~$j\in \{1,2,\cdots, m\}$~and all~$n$,
\begin{equation}
\sum_{i=1}^{m}Q_{n}(i,j) > 0,\quad \sum_{i=1}^{m}P_{n}(i,j) > 0 \quad \mathbb{P}\mbox{-almost surely}.
\end{equation}
C4 With positive $\mathbb{P}$-probability, the layer 0 is in one communication class.

\

\noindent{\it Known results.} Let us first  review some known results about the random walk in a random environment on the strip.

\noindent 1.{\it recurrence and transience.} If Condition C is satisfied, Theorem 1 in \cite{[BG00]} proved
 ~$\zeta_{n},~n\in \mathbb{Z}$~ of $m\times m$ matrices is the unique sequence
of stochastic matrices which satisfies the following system of equations:
\begin{equation}\label{ep}
\zeta_{n}=(I-Q_{n}\zeta_{n-1}-R_{n})^{-1}P_{n},\quad P-a.s.~n\in \mathbb{Z},
\end{equation}
and the enlarged sequence $(P_n,Q_n,R_n, \zeta_{n}), -\infty < n < \infty,$ is stationary and ergodic.

 ~Let
 \begin{equation}\label{p}
A_{n}:=(I-Q_{n}\zeta_{n-1}-R_{n})^{-1}Q_{n}\quad\mbox{and}\quad u_{n}:=(I-Q_{n}\zeta_{n-1}-R_{n})^{-1}\mathbf{1}
\end{equation}
and
\begin{equation}\label{lam}
\lambda^{+}:=\lim_{n\rightarrow \infty}\frac{1}{n}\log\parallel A_{n}A_{n-1}\cdots A_{1} \parallel,
\end{equation}
Theorem 2 in \cite{[BG00]} gave the criterion of recurrent and transient behavior for~$X_{n}=(\xi_{n},Y_{n})$. One of the cases is
\begin{equation}\label{rt}
\lim_{t\rightarrow \infty} \xi(t)= \infty, ~~\mathbb{P}-a.e.\quad \mbox{if and only if} \quad\lambda ^{+}< 0.
\end{equation}

\noindent 2.{\it exit probability.} Let~$\eta_{n}(i,j)$~be the probability of a random walk starting in~$(n,i)$~reaches the layer~$n+1$~at point~$(n,j)$~finally,~we usually called it the exiting probability. If the random walk is transient to the right, we have~$\eta_{n}=\zeta_{n},~P-a.e.$(see \cite{[Gol08]}, (1.15) ). And if Condition C is satisfied, ~$\zeta_{n}>0$~for~$P-a.s.~\omega$.

We only concentrate on random walks which are transient to the right in our paper.

\noindent 3.{\it  stationary
sequence of probability vectors $y_n$.}

If Condition C is satisfied then following limit exists for $P-a.s.$ $\omega$  (Lemma 1, \cite{[Gol08]}):
\begin{equation}\label{yn2}
\mathbf{y}_{n}:=\lim_{a\rightarrow -\infty}\mathbf{u}_{a}\zeta_{a}(\omega)\zeta_{a+1}(\omega)
\cdots\zeta_{n}(\omega).
\end{equation}
where~$\mathbf{u}_{a}$~is any sequence of row-vectors with non-negative components~$u_{a}(i)$, and
$\sum_{i=1}^{d}u_{a}(i)=1$. Note that the sequence~$\{\mathbf{y}_{n}\}$~is the unique solution of~$\mathbf{y}_{n}=\mathbf{y}_{n-1}\zeta_{n}$ in the class of probability vectors and it has the property
$y_n > 0$, which is a probability measure on~$\mathscr{D}=\{1,2,\cdots, d\}$~whose support is the whole~$\mathscr{D}$. It is clear that vectors $\mathbf{y}_n := \mathbf{y}(\omega_{\leq n}) $ form a stationary sequence.

\subsection{Statement of main results.}

We assume the walk ~$X_{n}=(\xi_{n},Y_{n})$ starts from layer $0$, the initial distribution~$P_{\omega}^{\mu}(\xi_{0}=0,~Y_{0}=i)=\mu_{\omega}(i), ~P-a.s. \omega,$~for any~$i\in \mathscr{D}$~with $\mu_{\omega}\in \mathcal {M}_{d}$.
In what follows, suppose Condition C is satisfied and~$\lambda ^{+}< 0$, i.e., we concentrate on random walks~$X_{n}=(\xi_{n},Y_{n})$~  transient to the right ~$X_{n}\rightarrow +\infty,~\mathbb{P}-a.s.$. ~In this case, suppose~$T_{0}=0$~and we have~$T_{k}<\infty,~\mathbb{P}-a.s.$~for any positive integer~$k\geq 1$. The aim of this paper is to calculate the hitting time ~$T_{1}=\inf\{i:~\xi(i)=1\}$~accurately in terms of the intrinsic branching structure within the walk. For~$n\leq1$,~define


$\mathbf{U}_{n}=(U_{n}^{1},U_{n}^{2},\cdots,U_{n}^{d})$,~where~$U_{n}^{i}~(1\leq i\leq d)$~is the number of steps from layer ~$n$~ to  layer ~$n-1$~at the site~$(n-1,i)$~before time~$T_{1}$.

$\mathbf{Z}_{n}=(Z_{n}^{1},Z_{n}^{2},\cdots,Z_{n}^{d})$,~where~$Z_{n}^{i}~(1\leq i\leq d)$~is the number of steps from layer ~$n$~ to the same layer at the site~$(n,i)$~before time~$T_{1}$.

And
\begin{equation}\label{duz}
|\mathbf{U}_{n}|=\sum_{i=1}^{d}U_{n}^{i}=\mathbf{U}_{n}\mathbf{1}
\quad \mbox{and} \quad |\mathbf{Z}_{n}|=\sum_{i=1}^{d}Z_{n}^{i}=\mathbf{Z}_{n}\mathbf{1}.
\end{equation}


All steps before $T_1$ can be recorded by~$\mathbf{U}_{n}$~and~$\mathbf{Z}_{n}$.~Since~$X_{n}\rightarrow +\infty,~\mathbb{P}-a.s.$,~if the random walk takes a step to the left from any layer~$n~(n\leq 0)$, it must come back finally from layer ~$n-1$~~to layer ~$n$~, so
$$T_{1}= 1+\sum_{n\leq 0}(2|\mathbf{U}_{n}|+|\mathbf{Z}_{n}|),$$
and the following theorem tells us that ~$\{|\mathbf{U}_{n}|, |\mathbf{Z}_{n}|, ~n\leq 1 \}$~is an inhomogeneous branching process with immigration. The exit probability $\eta_{n}$ plays an important role, when $X_{n}\rightarrow +\infty$, ~$\eta_{n}=\zeta_{n},~P-a.e.$, (see \cite{[Gol08]}, (1.15) ), which is given by (\ref{ep}).

\begin{theorem}\label{thm1}
 Assume Condition C is satisfied and ~$X_{n}\rightarrow +\infty,~\mathbb{P}-a.s.$, the initial distribution~$P_{\omega}^{\mu}(\xi_{0}=0,~Y_{0}=i)=\mu_{\omega}(i), ~P-a.s.. $ ~Then

\noindent (1) for~$P-a.s.~\omega$,~$\{|\mathbf{U}_{n}|,~n\leq 1\}$~and
~$\{|\mathbf{Z}_{n}|,~n\in\mathbb{Z}\}$~are  inhomogeneous branching processes with  immigration. The offspring distribution  is given by for $n\leq 0$
\begin{equation}\label{od1}
P_{\omega}^{\mu}\Big(|\mathbf{U}_{n}|=m \Big| \mathbf{U}_{n+1}=\mathbf{e}_{i}\Big )=\mathbf{e}_{i} [(I-R_{n})^{-1}  Q_{n}\zeta_{n-1}]^{m}(I-R_{n})^{-1}P_{n}\mathbf{1},
\end{equation}
\begin{equation}\label{od2}
P_{\omega}^{\mu}\Big(|\mathbf{Z}_{n}|=K \Big| \mathbf{U}_{n+1}=\mathbf{e}_{i}\Big)=\mathbf{e}_{i}[(I-Q_{n}\zeta_{n-1})^{-1} R_{n}]^{K} (I-Q_{n}\zeta_{n-1})^{-1} P_{n} \mathbf{1},
\end{equation}
with immigration
\begin{equation}\label{od3}
P_{\omega}^{\mu}\Big(\mathbf{U}_{1}=\mathbf{e}_{i} \Big )=\mu_{\omega}(i),  \ \ \ \  i\in \mathscr{D},
\end{equation}
where $\zeta_{n}=\eta_{n}$ (see \cite{[Gol08]}, (1.15) ) is exit probability, which is given by (\ref{ep}).

\noindent(2) The first hitting time $T_1$ is given by
\begin{equation}\label{t1t1}
T_{1}= 1+\sum_{n\leq 0}(2|\mathbf{U}_{n}|+|\mathbf{Z}_{n}|).
\end{equation}
\qed
\end{theorem}
\noindent{\textbf{Remark}}
(1) In Theorem \ref{thm1}, we restrict ourselves only to the trajectory of the walk $X_t$ for $t\in [0, T_1]$, and all the steps have been counted in $\{|\mathbf{U}_{n}|, |\mathbf{Z}_{n}|, ~n\leq 1 \}$ which formulate  a branching structure as (\ref{od1}) and  (\ref{od2}) with immigration  (\ref{od3}). After that, the trajectory of the walk $X_t$ follows the same structure. For example,  the trajectory of the walk $X_t$ for $t\in [T_1, T_2]$,  all the steps have been counted in $\{|\mathbf{U}_{n}|, |\mathbf{Z}_{n}|, ~n\leq 2 \}$ which formulate  a branching structure as (\ref{od1}) and  (\ref{od2}) with immigration  $P_{\omega}^{\mu}\Big(\mathbf{U}_{2}=\mathbf{e}_{i} \Big )=Y_{T_1}(i)$, and so on.

\noindent (2) Note that   it is  ``unsymmetrical" in  the branching structure (\ref{od1}) and  (\ref{od2}) between the ``father " and ``children". It can be explained as that we focus on the number of the ``children" but the individual of the ``father " (determine the probability).  \qed

As an immediate application of the branching structure, we can calculate the mean of the hitting time explicitly.
\begin{theorem}\label{thm2}
Assume Condition C is satisfied and ~$X_{n}\rightarrow +\infty,~\mathbb{P}-a.s.$, and the initial distribution~$P_{\omega}^{\mu}(\xi_{0}=0,~Y_{0}=i)=\mu_{\omega}(i), ~P-a.s.. $  Then
\begin{eqnarray}
&E{T_1}=&E({\overrightarrow{\mu_\omega}}(u_{0}+A_{0}u_{-1}+ \cdots +A_{0}A_{-1}\cdots A_{-k}u_{-k-1}+\cdots))\no,
\end{eqnarray}
where $A_n, u_n$ is given in  (\ref{p}).
\qed
\end{theorem}

Another application of the branching structure is to specify the density of   the absolutely continuous invariant measure for the ``environments viewed from the particle". Let us review the process discussed in Section 4 of $\cite{[Roi08]}$. Let ~$\overline{\omega_{n}}=\theta^{\xi_{n}}w$,~for~$n\geq 0$,~and consider the process~$Z_{n}:=(\overline{\omega_{n}},Y_{n})$, defined in~$(\Omega \times \mathscr{D},\mathscr{F}\otimes \mathscr{H})$,~where~$\mathscr{H}$~as the set of all subsets of~$\mathscr{D}$, and the initial distribution~$P_{\omega}^{\mu}(\xi_{0}=0,~Y_{0}=i)=\mu_{\omega}(i)=\mathbf{y}_{-1} (i) $ given by (\ref{yn2}). ~$(Z_{n})_{n\geq 0}$~is a Markov chain under~$\mathbb{P}^{\mu}$~with transition kernel
\begin{equation}
K(\omega,i;B,j)=P_{0}(i,j)I_{B}(\theta w)+R_{0}(i,j)I_{B}(w)+Q_{0}(i,j)I_{B}(\theta^{-1} w).
\end{equation}

\noindent Usually,~$Z_{n}=(\overline{\omega_{n}},Y_{n})$~be called as auxiliary Markov chain.

Let~$v_{p}=\frac{1}{\mathbb{E}T_{1}}$,~whenever~$\mathbb{E}T_{1} < \infty$.~For~$B \in \mathscr{F}$,~$i\in \mathscr{D}$,~define a probability measure~$Q$~on~$(\Omega \times \mathscr{D},\mathscr{F}\otimes \mathscr{H})$:
\begin{equation}\label{qb}
Q(B,i):=v_{p}\mathbb{E}\left(\sum_{n=0}^{T_{1}-1}
I_{B}(\overline{\omega_{n}})I_{Y_{n}}(i)\right).
\end{equation}
~$Q(\cdot)$~is a invariant measure under the Markov kernel~$K$ (Proposition 4.1, $\cite{[Roi08]}$).

Define a probability measure~$\overline{Q}(\cdot)$~on~$(\Omega,\mathscr{F})$~by setting
\begin{equation}\label{qh}
\overline{Q}(B):=Q(B,\mathscr{D}),\quad B \in \mathscr{F}.
\end{equation}
and let $~Q_{i}(B):=Q(B,i)~\mbox{for}~B \in \mathscr{F}$. Then both ~$Q_{i}(\cdot)$~ and ~$\overline{Q}(\cdot)$ are absolutely continuous with regard to~$P$ (Proposition 4.1, $\cite{[Roi08]}$), but where only the up bound of the  density have been proved. The branching structure enable us to specify the density completely in the following theorem.

\begin{theorem}\label{thm3}
Assume Condition C is satisfied and ~$X_{n}\rightarrow +\infty,~\mathbb{P}-a.s.$,   the initial distribution~$P_{\omega}^{\mu}(\xi_{0}=0,~Y_{0}=i)=\mu_{\omega}(i)=\mathbf{y}_{-1} (i) $ , ~for~${P}-a.s.~\omega$, ~and assume in addition that~$v_{p}>0$.~Then~$Q_{i}(\cdot)$~is absolutely continuous with regard to~$P$,~and~so is ~$\overline{Q}(\cdot)$. The density is given by
\begin{equation}\label{lwi}
\frac{dQ_{i}}{dP}=\Lambda_{\omega}^{(i)},
\end{equation}
\noindent where
\begin{equation}\label{lwi1}
\Lambda_{\omega}^{(i)}
 =v_{p} [\mathbf{\mu}_{\omega}\left(\widetilde{u}_{0}+\zeta_{0}A_{1}\widetilde{u}_{0}+ \zeta_{0}\zeta_{1}A_{2}A_{1}\widetilde{u}_{0}+\cdots \right)](i).
\end{equation}
and
\begin{equation}\label{lw1}
\frac{d\overline{Q}}{dP}=\Lambda_{\omega},
\end{equation}
\noindent where
\begin{equation}\label{lw2}
\Lambda_{\omega}
 = v_{p} [\mu_{\omega}\left(\widetilde{u}_{0}+\zeta_{0}A_{1}\widetilde{u}_{0}+ \zeta_{0}\zeta_{1}A_{2}A_{1}\widetilde{u}_{0}+\cdots \right)]\mathbf{1},
\end{equation}
where $\widetilde{u}_{n}:=(I-Q_{n}\zeta_{n-1}-R_{n})^{-1}$. \qed
\end{theorem}
\noindent{\textbf{Remark}}
(1) The first part of the Theorem  $\ref{thm3}$ is obtained in Proposition 4.1 of $\cite{[Roi08]}$. We will focus on the ``density" part only.

\noindent(2) As a by product, we can prove the LLN from two different method as the situation for the nearest RWRE on the line ($\cite{[Zei04]}$). On the one hand, If   $\overrightarrow{\mu_\omega}=\mathbf{y}_{-1} $, then $\{\tau_i: i\in N\}$ in (\ref{ht}) is a stationary and ergodic sequence variables (Lemma 3.2, $\cite{[Roi08]}$ ), so the LLN can be obtained from the hitting time decomposition; on the other hand, with the ``density" in hand, it is easy to obtain the LLN again from the point of view ``environments viewed from the particle". We omit the details of the proof.

\begin{corollary}\label{cor1}
Assume Condition C is satisfied and ~$X_{n}\rightarrow +\infty,~\mathbb{P}-a.s.$,   the initial distribution~$P_{\omega}^{\mu}(\xi_{0}=0,~Y_{0}=i)=\mu_{\omega}(i)=\mathbf{y}_{-1} (i) $ , ~for~${P}-a.s.~\omega$, ~and assume in addition that~$v_{p}>0$. Then $\mathbb{P}-a.s.$,
\begin{eqnarray}\label{xlln}
&&\lim_{n\rightarrow \infty} \frac{\xi(n)}{n}=\frac{1}{E\Big(\mathbf{y}_{-1}(u_{0}+A_{0}u_{-1}+ \cdots +A_{0}A_{-1}\cdots A_{-k}u_{-k-1}+\cdots)\Big)}.
\end{eqnarray}
\qed
\end{corollary}


\section{ Proofs  }

\setcounter{equation}{0}
\subsection{Intrinsic branching structure---{\it Proof of Theorem\ref{thm1}}.}
Assume Condition C is satisfied and ~$X_{n}\rightarrow +\infty,~\mathbb{P}-a.s.$, the initial distribution~$P_{\omega}^{\mu}(\xi_{0}=0,~Y_{0}=i)=\mu_{\omega}(i), ~P-a.s.. $  Note that~$T_{k}<\infty,~\mathbb{P}-a.s.$~for any positive integer~$k\geq 1$. We will analyze the trajectory of the walk, and restrict to  the first excursion between lay 0 to lay 1, i.e., the path of $X_k$ for $k\in [0,T_1]$. Define
 for~$n\leq 0$,
\begin{eqnarray}
\alpha_{n,0}&=&\min\{k\leq T_{1}:~X_{k}\in L_{n} \},\no\\
\beta_{n,0}&=&\min\{\alpha_{n,0}< k \leq T_{1}:~X_{k-1}\in L_{n},~X_{k}\in L_{n-1}\}.\no
\end{eqnarray}

\noindent And for~$b\geq 1$,
\begin{eqnarray}
\alpha_{n,b}&=&\min\{\beta_{n,b-1}< k\leq T_{1}:~X_{k}\in L_{n} \},\no\\
\beta_{n,b}&=&\min\{\alpha_{n,b}< k \leq T_{1}:~X_{k-1}\in L_{n},~X_{k}\in L_{n-1}\}.\no
\end{eqnarray}
(with the usual convention that the minimum over an empty set is~$+\infty$).

We refer to the time interval~$[\beta_{n,b-1},~\alpha_{n,b}]$~as the $b$-th excursion from~$n-1$~layer to~$n$~layer.

For any~$b\geq 0$,~any~$n\leq 0$,~and~$i\in \{1,2,\cdots,d\}$,~define
\begin{eqnarray}\label{uz}
U_{n,b}^{i}&:=&\sharp\{ k\geq 0:~X_{k-1}\in L_{n},~X_{k}=(n-1,i),~\beta_{n+1,b}<k<\alpha_{n+1,b+1}\},\\
Z_{n,b}^{i}&:=&\sharp\{ k\geq 0:~X_{k-1}\in L_{n},~X_{k}=(n,i),~\beta_{n+1,b}<k<\alpha_{n+1,b+1}\}.
\end{eqnarray}

Note that~$U_{n,b}^{i}$~is the number of steps from layer ~$n$~ to~$(n-1,i)$~during the~$b+1$-th excursion from layer ~$n$~ to layer ~$n+1$~, whereas ~$Z_{n,b}^{i}$~is the number of steps from layer ~$n$~ to~$(n,i)$~during the same excursion.

Define for~$n\leq0$~and~$i\in \{1,2,\cdots,d\}$,~$U_{n}^{i}:=\sum_{b\geq 0} U_{n,b}^{i}$,~then~$U_{n}^{i}$~is the number of steps from layer ~$n$~ to~$(n-1,i)$~before time~$T_{1}$.~Similarly define~$Z_{n}^{i}:=\sum_{b\geq 0}Z_{n,b}^{i}$. ~$\mathbf{U}_{n}=(U_{n}^{1},U_{n}^{2},\cdots,U_{n}^{d})$,
~and~$|\mathbf{U}_{n}|=\sum_{i=1}^{d}U_{n}^{i}=\mathbf{U}_{n}\mathbf{1}$;
~$\mathbf{Z}_{n}=(Z_{n}^{1},Z_{n}^{2},\cdots,Z_{n}^{d})$,
~and~$|\mathbf{Z}_{n}|=\sum_{i=1}^{d}Z_{n}^{i}=\mathbf{Z}_{n}\mathbf{1}$~which have been defined in (\ref{duz}).

By Markov property, we obtain
\begin{eqnarray}\label{bm11}
&& P_{\omega}^{\mu}\Big(|\mathbf{U}_{n}|=m, |\mathbf{Z}_{n}|=K \Big| \mathbf{U}_{n+1}=\mathbf{e_{i}}\Big)\no\\
&=& \mathbf{e}_{i} \sum\nolimits_{k_{0}+k_{1}+\cdots+k_{m}=K} R_{n}^{k_{0}} Q_{n}\zeta_{n-1} R_{n}^{k_{1}} \cdots Q_{n}\zeta_{n-1} R_{n}^{k_{m}} P_{n} \mathbf{1}.
\end{eqnarray}
where ~$\zeta_{n}=\eta_{n}$ is the exiting probability matrix (see ($\ref{ep}$)).

 In (\ref{bm11}), the path of an excursion is considered: the particle  start from  layer $n$ (given by $ \mathbf{U}_{n+1}=\mathbf{e_{i}}$), moves at layer $n$ by $|\mathbf{Z}_{n}|=K$ steps (each step with probability $R_{n}   $) and $|\mathbf{U}_{n}|=m$ steps from layer $n$ to layer $n-1$ ( but in the trajectory point, each ``down" step with probability $Q_{n}$ must connect with a path ``from layer $n-1$ finally goes back to layer $n$ " with probability $\zeta_{n-1}$), the last step of the excursion is from layer $n$ to layer $n+1$ with probability $P_{n}$.

The idea of (\ref{bm11}) is that we only care the number of the ``children", which lead to the ``unsymmetrical". Note that only the ``$ \mathbf{U}$ " type particles produce ``children". With a similar consideration, the branching mechanism can also be expressed as
\begin{eqnarray}\label{bm12}
&& P_{\omega}^{\mu}\Big(|\mathbf{U}_{n}|=m, |\mathbf{Z}_{n}|=K \Big| \mathbf{U}_{n+1}=\mathbf{e_{i}}\Big)\no\\
&=& \mathbf{e}_{i} \sum\nolimits_{m_{0}+m_{1}+\cdots+m_{K}=m} Q_{n}\zeta_{n-1}^{m_{0}} R_{n} Q_{n}\zeta_{n-1}^{m_{1}} \cdots R_{n} Q_{n}\zeta_{n-1}^{m_{K}} P_{n} \mathbf{1}.
\end{eqnarray}
In what follows, we will derive the marginal distribution of~$|\mathbf{U}_{n}|$~and~$|\mathbf{Z}_{n}|$~respectively. Let's discuss the marginal distribution of~$|\mathbf{U}_{n}|$ first, summarize over $K$ in (\ref{bm11}),

\begin{eqnarray}\label{mdu1}
&&P_{\omega}^{\mu}\Big(|\mathbf{U}_{n}|=m \Big| \mathbf{U}_{n+1}= \mathbf{e_{i}}\Big)\no\\
&=&\sum_{K=0}^{+\infty}P_{w}^{\mu}\Big(|\mathbf{U}_{n}|=m, |\mathbf{Z}_{n}|=K \Big| \mathbf{U}_{n+1}= \mathbf{e_{i}}\Big)\no\\
&=&\mathbf{e_{i}} \Big[\sum_{K=0}^{+\infty}  \sum\nolimits_{k_{0}+k_{1}+\cdots+k_{m}=K} R_{n}^{k_{0}} Q_{n}\zeta_{n-1} R_{n}^{k_{1}} \cdots Q_{n}\zeta_{n-1} R_{n}^{k_{m}} P_{n} \mathbf{1} \Big].
\end{eqnarray}
It's not hard to see
\begin{eqnarray}\label{mdu2}
&&\sum_{K=0}^{+\infty}  \sum\nolimits_{k_{0}+k_{1}+\cdots+k_{m}=K} R_{n}^{k_{0}} Q_{n}\zeta_{n-1} R_{n}^{k_{1}} \cdots Q_{n}\zeta_{n-1} R_{n}^{k_{m}}\no\\
&=&(I-R_{n})^{-1}  Q_{n}\zeta_{n-1}(I-R_{n})^{-1}  \cdots Q_{n}\zeta_{n-1}(I-R_{n})^{-1} \no\\
&=&[(I-R_{n})^{-1}  Q_{n}\zeta_{n-1}]^{m}(I-R_{n})^{-1}.
\end{eqnarray}
Taking together (\ref{mdu1}) and (\ref{mdu2}), derives the marginal distribution of~$|\mathbf{U}_{n}|$,
\begin{equation}\label{mdu3}
P_{\omega}^{\mu}\Big(|\mathbf{U}_{n}|=m \Big| \mathbf{U}_{n+1}=\mathbf{e_{i}}\Big)=\mathbf{e_{i}} [(I-R_{n})^{-1}  Q_{n}\zeta_{n-1}]^{m} (I-R_{n})^{-1}P_{n} \mathbf{1}.
\end{equation}
For the marginal distribution of~$|\mathbf{Z}_{n}|$, summarize over $m$ in (\ref{bm12}), we have
\begin{eqnarray}\label{mdz}
&& P_{\omega}^{\mu}\Big(|\mathbf{Z}_{n}|=K \Big| \mathbf{U}_{n+1}=\mathbf{e_{i}}\Big)\no\\
&=&\mathbf{e_{i}} \sum_{m=0}^{+\infty}\sum\nolimits_{m_{0}+m_{1}+\cdots+m_{K}=m} Q_{n}\zeta_{n-1}^{m_{0}} R_{n} Q_{n}\zeta_{n-1}^{m_{1}} \cdots R_{n} Q_{n}\zeta_{n-1}^{m_{K}} P_{n} \mathbf{1}\no\\
&=&\mathbf{e_{i}} (I-Q_{n}\zeta_{n-1})^{-1} R_{n} (I-Q_{n}\zeta_{n-1})^{-1} \cdots R_{n}(I-Q_{n}\zeta_{n-1})^{-1}  P_{n} \mathbf{1}\no\\
&=&\mathbf{e_{i}} [(I-Q_{n}\zeta_{n-1})^{-1} R_{n}]^{K} (I-Q_{n}\zeta_{n-1})^{-1} P_{n} \mathbf{1}.
\end{eqnarray}
 Complete the proof of part (1) of  Theorem (\ref{thm1}); and part (2) is immediate. \qed

\noindent{\textbf{Remark.}} \quad
From the marginal distribution, we also can test of the validity of the branching structure. In fact,
~\beqnn
  \sum_{m=0}^{+\infty} P_{\omega}^{\mu}\Big(|\mathbf{U}_{n}|=m \Big| \mathbf{U}_{n+1}=\mathbf{e}_{i}\Big)
 & = & \mathbf{e_{i}}[ \sum_{m=0}^{+\infty}  [(I-R_{n})  Q_{n}\zeta_{n-1}]^{m} ](I-R_{n})^{-1}P_{n} \mathbf{1}\\
 & = & \mathbf{e_{i}}[ I-(I-R_{n})  Q_{n}\zeta_{n-1}]^{-1} ](I-R_{n})^{-1}P_{n} \mathbf{1}\\
 & = & \mathbf{e_{i}}[ (I-R_{n}) - Q_{n}\zeta_{n-1}]^{-1} ]P_{n} \mathbf{1}\\
 & = & \mathbf{e_{i}}\zeta_{n} \mathbf{1}=1.
  \eeqnn

\subsection{{\bf $E{T_1}$}---{\it Proof of Theorem \ref{thm2}}}

The random walk~$X_{n}=(\xi_{n},Y_{n})$~starts from layer ~$0$~ with the initial distribution~$\mathbf{\mu}_{\omega}$. With the branching structure in hand, we can calculate the mean of the first hitting time $T_1$. We discuss it by four steps as follows.

\noindent{\textbf{Step 1.}} $E_{\omega}^{\mu}\Big(| \mathbf{U}_{n}| \Big|  \mathbf{U}_{n+1}=\mathbf{e_{i}}\Big)$ and $E_{\omega}^{\mu}\Big(| \mathbf{Z}_{n}| \Big|  \mathbf{U}_{n+1}=\mathbf{e_{i}}\Big)$.

\noindent From (\ref{od1}) of Theorem \ref{thm1},
\begin{eqnarray}\label{expu}
 E_{\omega}^{\mu}\Big(| \mathbf{U}_{n}| \Big|  \mathbf{U}_{n+1}=\mathbf{e_{i}}\Big)
 &=& \sum_{m=0}^{+\infty} m P_{\omega}^{\mu}\Big(|\mathbf{U}_{n}|=m \Big| \mathbf{U}_{n+1}=\mathbf{e_{i}}\Big)\no\\
&=&\mathbf{e_{i}} \sum_{m=1}^{+\infty} m [(I-R_{n})^{-1} Q_{n}\zeta_{n-1}]^{m} (I-R_{n})^{-1}P_{n} \mathbf{1}.
\end{eqnarray}

\noindent To process the calculation, we need the following
\begin{lemma}\label{mat}
For matrix $B$,~$I-B$~is non-degenerate, then~$\sum_{m=1}^{+\infty} m B^{m}=B(I-B)^{-2}$.
\end{lemma}

\proof $$\sum_{m=1}^{+\infty} m B^{m}=(B+2B^{2}+3B^{3}+\cdots)=B(I+2B+3B^{2}+\cdots),$$
and
\begin{equation*}
\begin{split}
(I-B)^{-2}=((I-B)^{-1})^{2}=(\sum_{m=1}^{+\infty}B^{n})^{2}
=(\sum_{m=1}^{+\infty}B^{n})(\sum_{m=1}^{+\infty}B^{n})
= (I+2B+3B^{2}+4B^{3}\cdots).
\end{split}
\end{equation*}
Thus~$\sum_{m=1}^{+\infty} m B^{m}=B(I-B)^{-2}$. \qed

\

 Let~$B=(I-R_{n})^{-1} Q_{n}\zeta_{n-1}$, (\ref{expu}) can be continued as
\begin{eqnarray}\label{expu2}
E_{\omega}^{\mu}\Big(| \mathbf{U}_{n}| \Big|  \mathbf{U}_{n+1}=\mathbf{e_{i}}\Big) &=&\mathbf{e_{i}}  (I-R_{n})^{-1} Q_{n}\zeta_{n-1} [I-(I-R_{n})^{-1} Q_{n}\zeta_{n-1}]^{-2} (I-R_{n})^{-1}P_{n} \mathbf{1} \no \\
 &=&\mathbf{e_{i}}  (I-Q_{n}\zeta_{n-1}-R_{n})^{-1}Q_{n}\zeta_{n-1}\zeta_{n} \mathbf{1} \\
 &=&\mathbf{e_{i}}  A_{n} \mathbf{1}.
\end{eqnarray}
The second equality (\ref{expu2}) need a series calculations about the matrix which we leave it as Appendix, where $A_n$ is given in (\ref{p}). Similarly,
\begin{eqnarray}\label{expz}
E_{\omega}^{\mu}\Big(| \mathbf{Z}_{n}| \Big|  \mathbf{U}_{n+1}=\mathbf{e_{i}}\Big)
&=&\sum_{K=0}^{+\infty} \mathbf{e_{i}}K [(I-Q_{n}\zeta_{n-1})^{-1}R_{n} ]^{K} (I-Q_{n}\zeta_{n-1})^{-1}P_{n} \mathbf{1}\no\\
&=&\mathbf{e_{i}}  (I-Q_{n}\zeta_{n-1}-R_{n})^{-1}R_{n}\zeta_{n} \mathbf{1}\no\\
&=&\mathbf{e_{i}}  (I-Q_{n}\zeta_{n-1}-R_{n})^{-1}R_{n} \mathbf{1}.
\end{eqnarray}
As a result, we have
\begin{eqnarray}\label{exp}
&&E_{\omega}^{\mu}\Big(| \mathbf{U}_{n}| \Big| \mathbf{U}_{n+1}\Big)=\mathbf{U}_{n+1}A_{n} \mathbf{1},\\
&&E_{\omega}\Big(|\mathbf{Z}_{n}|  \Big| \mathbf{U}_{n+1}\Big)=\mathbf{U}_{n+1}(I-Q_{n}\zeta_{n-1}-R_{n})^{-1}R_{n} \mathbf{1}.
\end{eqnarray}

 \noindent{\textbf{Step 2.}} Steps visited on layer $n$.

For any~$n\leq 0$,~define
\begin{equation}\label{n1}
N_{n}^{i}=\sharp\{ k \in [0,T_{1}):~X_{k}=(n,i)\}.
\end{equation}
Note that~$N_{n}^{i}$~is the number of steps visited~$(n,i)$~before time~$T_{1}$. Let ~$\mathbf{N}_{n}=(N_{n}^{1},N_{n}^{2},
\cdots,N_{n}^{d})$~
and~$|\mathbf{N}_{n}|=\sum_{i=1}^{d}N_{n}^{i}=\mathbf{N}_{n}\mathbf{1}$.

Define a vector valued random variable~$\mathbf{U'}_{n}$~where~$\mathbf{U'}$$_{n}^{i},~1 \leq i \leq d$~is the number of steps from layer ~$n-1$~ to~$(n,i)$.~Then
\begin{equation}\label{n2}
|\mathbf{N}_{n}|=|\mathbf{U'}_{n}|+|\mathbf{Z}_{n}|+
|\mathbf{U}_{n+1}|,\quad \mathbb{P}-a.s..
\end{equation}

For another perspective,~$T_{1}=\sum_{n\leq 0}(|\mathbf{N}_{n}|), ~ \mathbb{P}-a.s..$
Since~$X_{n}\rightarrow +\infty,~\mathbb{P}-a.s.$,~if the random walk takes a step to the left from any layer~$n~(n\leq 0)$~to layer ~$n-1$~, it must come back finally from layer ~$n-1$~~to layer ~$n$~, therefore $|\mathbf{U}_{n}|=|\mathbf{U'}_{n}|,\quad \mathbb{P}-a.s..$

Together with (\ref{n2}), we have
\begin{eqnarray}
&& E_{\omega}^{\mu}(| \mathbf{N}_{n}|)
=E_{\omega}^{\mu}(|\mathbf{U}_{n}|+|\mathbf{Z}_{n}|+|\mathbf{U}_{n+1}|)\no\\
&=&E_{\omega}^{\mu}\Big[E_{\omega}^{\mu}(|\mathbf{U}_{n}|\Big|  \mathbf{U}_{n+1})+E_{\omega}^{\mu}(|\mathbf{Z}_{n}|\Big|  \mathbf{U}_{n+1})+E_{\omega}^{\mu}(|\mathbf{U}_{n+1}| \Big|  \mathbf{U}_{n+1})\Big]\no.
\end{eqnarray}

By using (\ref{exp}), one can calculate the quenched expectation of~$| \mathbf{N}_{n}|$~as
\begin{eqnarray}\label{en}
&& E_{\omega}^{\mu}(| \mathbf{N}_{n}|)
 =E_{\omega}^{\mu}[\mathbf{U}_{n+1}A_{n} \mathbf{1}+ \mathbf{U}_{n+1}(I-Q_{n}\zeta_{n-1}-R_{n})^{-1}R_{n} \mathbf{1}+\mathbf{U}_{n+1}\mathbf{1}]\no\\
&=&E_{\omega}^{\mu}[\mathbf{U}_{n+1}(I-Q_{n}\zeta_{n-1}-R_{n})^{-1}
 (Q_{n}\zeta_{n-1}+
 R_{n}+I-Q_{n}\zeta_{n-1}-R_{n})\mathbf{1}]\no\\
&=&E_{\omega}^{\mu}(\mathbf{U}_{n+1})(I-Q_{n}\zeta_{n-1}-R_{n})^{-1}\mathbf{1}.
\end{eqnarray}

 \noindent{\textbf{Step 3.}}
The next object is to discuss~$E_{\omega}^{\mu}(\mathbf{U}_{n+1})$.

Define a probability matrix~$B_{m}$,~where~$B_{m}(i,j)$~is the probability of a particle starting from~$(n+1,i)$,~takes more than~$m$~steps to the layer $n$, and the~$m$-th step located at~$(n,j)$. $B_{m}(i,j)$~can be expressed by our branching structure
\begin{eqnarray}\label{bm1}
B_{m}(i,j)
&=& \mathbf{e_{i}}\Big[\sum_{\overline{K}=0}^{+\infty}  \sum\nolimits_{k_{0}+k_{1}+\cdots+k_{m-1}=\overline{K}} R_{n+1}^{k_{0}} Q_{n+1}\zeta_{n} R_{n+1}^{k_{1}}
Q_{n+1}\zeta_{n} R_{n+1}^{k_{2}} \cdots Q_{n+1}\zeta_{n} R_{n+1}^{k_{m-1}} Q_{n+1}\Big]\mathbf{e_{j}}\no\\
&=& \mathbf{e_{i}} [(I- R_{n+1}) Q_{n+1}\zeta_{n}]^{m-1} (I- R_{n+1})^{-1} Q_{n+1}\mathbf{e_{j}}.
\end{eqnarray}
Let
$$
\widetilde{P}_{i,j}^{m}:={B}_{m}(i,j) - {B}_{m+1}(i,j),
$$
be the probability of a particle starts from~$(n+1,i)$, the~$m$-th step takes to the left and located at~$(n,j)$ . We have
\begin{eqnarray}\label{bm0}
E_{\omega}^{\mu}({U}_{n+1}^{j} \Big| \mathbf{U}_{n+2}=\mathbf{e_{i}})
&= & \sum_{m=1}^{+\infty} m \widetilde{P}_{i,j}^{m}
= \mathbf{e_{i}}\sum_{m=1}^{+\infty} m({B}_{m}-{B}_{m+1})\mathbf{e_{j}}\no\\
&= & \mathbf{e_{i}}\sum_{m=1}^{+\infty}B_{m}\mathbf{e_{j}}.
\end{eqnarray}
Combine with (\ref{bm1}),
\begin{eqnarray}\label{bm00}
E_{\omega}^{\mu}(\mathbf{U}_{n+1} | \mathbf{U}_{n+2})&=&\mathbf{U}_{n+2}   \sum_{m=1}^{+\infty} {B}_{m} \no \\
&=& \mathbf{U}_{n+2}   \sum_{m=1}^{+\infty} [(I- R_{n+1}) Q_{n+1}\zeta_{n}]^{m-1} (I- R_{n+1})^{-1} Q_{n+1} \no\\
&= &   \mathbf{U}_{n+2}  (I- Q_{n+1}\zeta_{n}-R_{n+1})^{-1} Q_{n+1} \no\\
&=& \mathbf{U}_{n+2}A_{n+1}.
\end{eqnarray}
Then
\begin{equation}
E_{\omega}^{\mu}(\mathbf{U}_{n+1})
= E_{\omega}^{\mu}[E_{\omega}^{\mu}(\mathbf{U}_{n+1}| \mathbf{U}_{n+2})]
= E_{\omega}^{\mu}(\mathbf{U}_{n+2})A_{n+1},
\end{equation}
where $A_n$ is given in (\ref{p}). By recursive argument, we obtain
\begin{eqnarray}\label{re}
E_{\omega}^{\mu}(\mathbf{U}_{n+1})
& = & E_{\omega}^{\mu}(\mathbf{U}_{n+3})A_{n+2}A_{n+1}\no\\
& = & \cdots\no\\
& = & E_{\omega}^{\mu}(\mathbf{U}_{1})A_{0}A_{-1}A_{-2} \cdots A_{n+2}A_{n+1}.
\end{eqnarray}

 \noindent{\textbf{Step 4.}}{\it Calculate $\mathbb{E}(T_{1})$}.  It follows from
(\ref{en}) and (\ref{re}) that,
\begin{eqnarray}\label{lt1}
 \mathbb{E}(T_{1})&=&E(E_{\omega}(T_{1}))
= E\left(\sum_{n\leq 0}E_{\omega}(|\mathbf{N}_{n}|)\right)\no\\
&=& E\left(\sum_{n\leq 0}E_{\omega}(\mathbf{U}_{n+1})(I-Q_{n}\zeta_{n-1}-R_{n})^{-1}\mathbf{1}\right)\no\\
&=&E\left(\sum_{n\leq 0}E_{\omega}(\mathbf{U}_{1})A_{0}A_{-1}A_{-2} \cdots A_{n+2}A_{n+1}(I-Q_{n}\zeta_{n-1}-R_{n})^{-1}\mathbf{1}\right)\no\\
&=&E(E_{\omega}^{\mu}(\mathbf{U}_{1})(u_{0}+A_{0}u_{-1}+ \cdots +A_{0}A_{-1}\cdots A_{-k}u_{-k-1}+\cdots))\no\\
&=&E(\overrightarrow{\mu_\omega}(u_{0}+A_{0}u_{-1}+ \cdots +A_{0}A_{-1}\cdots A_{-k}u_{-k-1}+\cdots)).
\end{eqnarray}
where $A_n, u_n$ is given in  (\ref{p}). \qed


\subsection{Density of the absolutely continuous invariant measure--{\it Proof of Theorem \ref{thm3}}}

Let us review the process discussed in Section 4 of $\cite{[Roi08]}$. From the point of view ``environments viewed from the particle", let ~$\overline{\omega_{n}}=\theta^{\xi_{n}}w$,~for~$n\geq 0$,~and consider the process~$Z_{n}:=(\overline{\omega_{n}},Y_{n})$, defined in~$(\Omega \times \mathscr{D},\mathscr{F}\otimes \mathscr{H})$,~where~$\mathscr{H}$~as the set of all subsets of~$\mathscr{D}$, and the initial distribution~$P_{\omega}^{\mu}(\xi_{0}=0,~Y_{0}=i)=\mu_{\omega}(i)=\mathbf{y}_{-1} (i) $ given by (\ref{yn2}). ~$(Z_{n})_{n\geq 0}$~is a Markov chain under~$\mathbb{P}^{\mu}$~with transition kernel
\begin{equation}
K(\omega,i;B,j)=P_{0}(i,j)I_{B}(\theta w)+R_{0}(i,j)I_{B}(w)+Q_{0}(i,j)I_{B}(\theta^{-1} w).
\end{equation}

Let~$v_{p}=\frac{1}{\mathbb{E}T_{1}}$,~whenever~$\mathbb{E}T_{1} < \infty$.~For~$B \in \mathscr{F}$,~$i\in \mathscr{D}$,~define a probability measure~$Q$~on~$(\Omega \times \mathscr{D},\mathscr{F}\otimes \mathscr{H})$:
\begin{equation}
Q(B,i):=v_{p}\mathbb{E}\left(\sum_{n=0}^{T_{1}-1}
I_{B}(\overline{\omega_{n}})I_{Y_{n}}(i)\right)=v_{p}\sum_{j\in \mathscr{D}}E_{p}\left(\mu_{\omega}(j)E_{\omega}^{j}
\left(\sum_{n=0}^{T_{1}-1}
I_{B}(\theta^{\xi_{n}}\omega)I_{Y_{n}}(i)\right)\right)\no
\end{equation}
~$Q(\cdot)$~is a invariant measure under the Markov kernel~$K$ (Proposition 4.1, $\cite{[Roi08]}$).

Define a probability measure~$\overline{Q}(\cdot)$~on~$(\Omega,\mathscr{F})$~by setting
\begin{equation}\label{qh1}
\overline{Q}(B):=Q(B,\mathscr{D}),\quad B \in \mathscr{F}.
\end{equation}
and let $~Q_{i}(B):=Q(B,i)~\mbox{for}~B \in \mathscr{F}$. Then both ~$Q_{i}(\cdot)$~ and ~$\overline{Q}(\cdot)$ are absolutely continuous with regard to~$P$ (Proposition 4.1, $\cite{[Roi08]}$), but where only the up bound of the  density have been proved.

For~$m\leq 0,~i\in \mathscr{D}$,~define~$N_{m}^{i}$~as in (\ref{n1}):
\begin{equation}
N_{m}^{i}:=\left\{\sharp n \in [0,T_{1}):~\xi_{n}=m,~Y_{n}=i \right\}.\no
\end{equation}
Note that for any bounded measurable function~$f:~\Omega \rightarrow \mathbb{R}$~and~$\forall~ i \in \mathscr{D}$,
\begin{eqnarray}\label{inme}
\int_{\Omega}f(\omega)Q(d\omega,i)
&=&v_{p}\sum_{n=0}^{+\infty}\mathbb{E}^{\mu}(f(\overline{w_{n}})
;~Y_{n}=i,~T_{1}>n)\no\\
&=& v_{p}\sum_{m\leq 0}\mathbb{E}^{\mu}(f(\theta^{m} \omega)N_{m}^{i})\no\\
&=& v_{p}E_{p}\left(\sum_{k\in \mathscr{D}}\mu_{\omega}(k) \sum_{m\leq 0}\left(f(\theta^{m} \omega) E_{\omega}^{k}N_{m}^{i}\right)\right)\no\\
&=&v_{p}E_{p}\left(f(\omega) \sum_{k\in \mathscr{D}} \sum_{m\leq 0} \mu_{\theta^{-m}\omega}(k) E_{\theta^{-m}\omega}^{k}N_{m}^{i}\right).
\end{eqnarray}

Therefore,~$Q_{i}$~is absolutely continuous with respect to~$P$,~and also~$\overline{Q}$~is absolutely continuous with respect to~$P$.~And the density is
\begin{eqnarray}\label{dens}
&&\Lambda_{\omega}^{(i)}:=\frac{dQ_{i}}{d P}
 =v_{p} \sum_{k\in \mathscr{D}} \sum_{m\leq 0} \mu_{\theta^{-m}\omega}(k) E_{\theta^{-m}\omega}^{k}(N_{m}^{i}) \\
&&\Lambda_{\omega}:=\frac{d\overline{Q}}{d P}
 = v_{p} \sum_{k\in \mathscr{D}} \sum_{m\leq 0} \mu_{\theta^{-m}\omega}(k) E_{\theta^{-m}\omega}^{k}(\mathbf{N}_{m}\mathbf{1}).
\end{eqnarray}

We intend to spesity the density~$\Lambda_{\omega}^{(i)}$~and~$\Lambda_{\omega}$~by branching structure. Note that $\mu_{\omega}(i)=\mathbf{y}_{-1} (i) $ given by (\ref{yn2}),  and $\zeta_{-n}=\eta_{-n}$ is the exit probability,
\begin{eqnarray}
 \mu_{\theta \omega}
&=& \lim _{n\rightarrow\infty}\mathbf{e}_{i}\zeta_{-n}(\theta \omega)\cdots \zeta_{-2}(\theta \omega)\zeta_{-1}(\theta \omega)\no\\
&=& \lim _{n\rightarrow\infty}\mathbf{e}_{i}\zeta_{-n+1}(\omega)\cdots \zeta_{-1}(\omega)\zeta_{0}(\omega)\no\\
&=& \mu_{\omega}\zeta_{0}(\omega)\no.
\end{eqnarray}

\noindent Thus
\begin{equation}
\sum_{k\in \mathscr{D}} \mu_{\theta \omega}(k) E_{\theta \omega}^{k}(N_{-1}^{i})
 =\sum_{k\in \mathscr{D}} \mu_{\theta \omega}(k) E_{\omega}^{k}(N_{0}^{i})
 =\sum_{k\in \mathscr{D}} \mu_{\omega}\zeta_{0}(k) E_{\omega}^{k}(N_{0}^{i}).\no
\end{equation}

\noindent Similarly, for~$m\leq 0$~and~$ i\in \mathscr{D}$
\begin{equation}\label{soe}
\sum_{k\in \mathscr{D}} \mu_{\theta^{-m}\omega}(k) E_{\theta^{-m}\omega}^{k}(N_{m}^{i})
 =\sum_{k\in \mathscr{D}} \mu_{\omega}\zeta_{0}\zeta_{1} \cdots \zeta_{-m-1}(k) E_{\omega}^{k}(N_{0}^{i}).
\end{equation}

The following lemma is closely related to branching structure.
\begin{lemma}\label{lem3}
For~$n < 0$,
\begin{equation}
E_{\omega}(\mathbf{N}_{n})=\mu_{\omega}A_{0}A_{-1}\cdots A_{n+1}\widetilde{u}_{n}.
\end{equation}
\end{lemma}

\proof
Due to the definition of~$\mathbf{N}_{n}$,~$\mathbf{U}_{n}'$,~$\mathbf{Z}_{n}$~and~$\mathbf{U}_{n+1}$,
\begin{equation}
\mathbf{N}_{n}=\mathbf{U}_{n}'+\mathbf{Z}_{n}+\mathbf{U}_{n+1}.\no
\end{equation}

\noindent Recall the branching structure and by similarly argument as in the proof of Theorem \ref{thm2}, we obtain that
\begin{eqnarray}
E_{\omega}(\mathbf{U}_{n}')
&=& \mathbf{U}_{n+1} \sum_{m=1}^{+\infty}\left((I-R_{n})^{-1}Q_{n}\zeta_{n-1}\right)^{m-1}
 (I-R_{n})^{-1}Q_{n}\zeta_{n-1}\no\\
&=& \mathbf{U}_{n+1} \left((I-Q_{n}\zeta_{n-1}-R_{n}\right)^{-1}Q_{n}\zeta_{n-1}
= \mathbf{U}_{n+1} A_{n}\zeta_{n-1}\no,
\end{eqnarray}
and
\begin{eqnarray}
E_{\omega}(\mathbf{Z}_{n})
&=&  \mathbf{U}_{n+1} \sum_{K=0}^{+\infty}\left((I-Q_{n}\zeta_{n-1})^{-1}R_{n}\right)^{K}
 (I-Q_{n}\zeta_{n-1})^{-1}R_{n}\no\\
&=&  \mathbf{U}_{n+1} \left(I-Q_{n}\zeta_{n-1}-R_{n}\right)^{-1}R_{n}\no.
\end{eqnarray}
Thus
\begin{eqnarray}
E_{\omega}(\mathbf{N}_{n}\mid \mathbf{U}_{n+1})
&=& E_{\omega}(\mathbf{N}_{n}\mid \mathbf{U}_{n+1})\no\\
&=& E_{\omega}(\mathbf{U}_{n}'\mid \mathbf{U}_{n+1})+E_{\omega}(\mathbf{Z}_{n}\mid \mathbf{U}_{n+1})+E_{\omega}(\mathbf{U}_{n+1}\mid \mathbf{U}_{n+1})\no\\
&=& E_{\omega}(\mathbf{U}_{n+1}\left[(I-Q_{n}\zeta_{n-1}-R_{n})^{-1}
 Q_{n}\zeta_{n-1}+(I-Q_{n}\zeta_{n-1}-R_{n})^{-1}R_{n}+I \right])\no\\
&=& E_{\omega}(\mathbf{U}_{n+1})(I-Q_{n}\zeta_{n-1}-R_{n})^{-1}\no.
\end{eqnarray}
Together with the fact
\begin{equation}
E_{\omega}(\mathbf{U}_{n+1})=\mathbf{U}_{n+2}A_{n+1},\no
\end{equation}
we have
\begin{eqnarray}
E_{\omega}(\mathbf{N}_{n})
&=&\mu_{\omega}A_{0}A_{-1}\cdots A_{n+1}(I-Q_{n}\zeta_{n-1}-R_{n})^{-1}\no\\
&=& \mu_{\omega}A_{0}A_{-1}\cdots A_{n+1}\widetilde{u}_{n}\no.
\end{eqnarray}
Then Lemma \ref{lem3} follows.\qed

It follows from equation (\ref{soe}) and lemma \ref{lem3} that
\begin{eqnarray}
\sum_{k\in \mathscr{D}} \mu_{\theta^{-m}\omega}(k) E_{\theta^{-m}\omega}^{k}(N_{m}^{i})
&=& \sum_{k\in \mathscr{D}} \mu_{\omega}\zeta_{0}\zeta_{1} \cdots \zeta_{-m-1}(k) E_{\omega}^{k}(N_{0}^{i})\no\\
&=& \mu_{\omega}\zeta_{0}\zeta_{1} \cdots \zeta_{-m-1}A_{-m}A_{-m-1}\cdots A_{2}A_{1}\widetilde{u}_{0}(i).
\end{eqnarray}

Thus
\begin{eqnarray}
\Lambda_{\omega}^{(i)}=\frac{dQ_{i}}{dP}
&=& v_{p} \sum_{k\in \mathscr{D}} \sum_{m\leq 0} \mu_{\theta^{-m}\omega}(k) E_{\theta^{-m}\omega}^{k}(N_{m}^{i})\no\\
&=& v_{p} \sum_{m\leq 0}[ \mu_{\omega}\zeta_{0}\zeta_{1} \cdots \zeta_{-m-1}A_{-m}A_{-m-1}\cdots A_{2}A_{1}\widetilde{u}_{0}](i)\no\\
&=& v_{p}[ \mu_{\omega}\left(\widetilde{u}_{0}+\zeta_{0}A_{1}\widetilde{u}_{0}+ \zeta_{0}\zeta_{1}A_{2}A_{1}\widetilde{u}_{0}+\cdots \right)](i)\no.
\end{eqnarray}
and similarly
\begin{eqnarray}
\frac{d\overline{Q}}{dP}=\Lambda_{\omega}
&=& v_{p} \sum_{k\in \mathscr{D}} \sum_{m\leq 0} \mu_{\theta^{-m}\omega}(k) E_{\theta^{-m}\omega}^{k}(\mathbf{N}_{m}\mathbf{1})\no\\
&=& v_{p} \sum_{m\leq 0} [\mu_{\omega}\zeta_{0}\zeta_{1} \cdots \zeta_{-m-1}A_{-m}A_{-m-1}\cdots A_{2}A_{1}\widetilde{u}_{0}]\mathbf{1} \no\\
&=& v_{p} [\mu_{\omega}\left(\widetilde{u}_{0}+\zeta_{0}A_{1}\widetilde{u}_{0}+ \zeta_{0}\zeta_{1}A_{2}A_{1}\widetilde{u}_{0}+\cdots \right)]\mathbf{1}\no.
\end{eqnarray}
\qed


{\center\section*{Appendix}}
\setcounter{equation}{0}
The following calculation is needed in (\ref{expu2}), it is a details calculations on the matrix.
~\beqnn
 & & E_{w}^{\mu}(| \mathbf{U}_{n}| \mid  \mathbf{U}_{n+1}=\mathbf{e_{i}})\\
 & = & \sum_{m=0}^{+\infty} m P_{w}^{\mu}(|\mathbf{U}_{n}|=m \mid \mathbf{U}_{n+1}=\mathbf{e_{i}})\\
 & = & \sum_{m=0}^{+\infty} \mathbf{e_{i}}m [(I-R_{n})^{-1} Q_{n}\zeta_{n-1}]^{m} (I-R_{n})^{-1}P_{n} \mathbf{1}\\
 & = & \mathbf{e_{i}} \sum_{m=1}^{+\infty} m [(I-R_{n})^{-1} Q_{n}\zeta_{n-1}]^{m} (I-R_{n})^{-1}P_{n} \mathbf{1}\\
  & = & \mathbf{e_{i}}  (I-R_{n})^{-1} Q_{n}\zeta_{n-1} [I-(I-R_{n})^{-1} Q_{n}\zeta_{n-1}]^{-2} (I-R_{n})^{-1}P_{n} \mathbf{1}\\
 & = & \mathbf{e_{i}}  (I-R_{n})^{-1} Q_{n}\zeta_{n-1} [I-(I-R_{n})^{-1} Q_{n}\zeta_{n-1}]^{-1} [I-(I-R_{n})^{-1} Q_{n}\zeta_{n-1}]^{-1} (I-R_{n})^{-1}P_{n} \mathbf{1}\\
 & = & \mathbf{e_{i}} \{[(I-R_{n})^{-1} Q_{n}\zeta_{n-1}]^{-1}\}^{-1} [I-(I-R_{n})^{-1} Q_{n}\zeta_{n-1}]^{-1} [I-(I-R_{n})^{-1} Q_{n}\zeta_{n-1}]^{-1} (I-R_{n})^{-1}P_{n}  \mathbf{1}\\
 & = & \mathbf{e_{i}} \{[I-(I-R_{n})^{-1} Q_{n}\zeta_{n-1}][(I-R_{n})^{-1} Q_{n}\zeta_{n-1}]^{-1}\}^{-1}[I-(I-R_{n})^{-1} Q_{n}\zeta_{n-1}]^{-1} (I-R_{n})^{-1}P_{n}  \mathbf{1}\\
 & = & \mathbf{e_{i}}  \{[(I-R_{n})^{-1} Q_{n}\zeta_{n-1}]^{-1}-I\}^{-1} [I-(I-R_{n})^{-1} Q_{n}\zeta_{n-1}]^{-1} (I-R_{n})^{-1}P_{n}\mathbf{1}\\
 & = & \mathbf{e_{i}} \{[I-(I-R_{n})^{-1} Q_{n}\zeta_{n-1}][[(I-R_{n})^{-1} Q_{n}\zeta_{n-1}]^{-1}-I]\}^{-1} (I-R_{n})^{-1}P_{n} \mathbf{1}\\
 & = & \mathbf{e_{i}}  \{(I-R_{n})[I-(I-R_{n})^{-1} Q_{n}\zeta_{n-1}][[(I-R_{n})^{-1} Q_{n}\zeta_{n-1}]^{-1}-I]\}^{-1}P_{n} \mathbf{1}\\
 & = & \mathbf{e_{i}} \{[(I-R_{n}) -Q_{n}\zeta_{n-1}][[(I-R_{n})^{-1} Q_{n}\zeta_{n-1}]^{-1}-I]\}^{-1}P_{n}  \mathbf{1}\\
 & = & \mathbf{e_{i}}  \{[(I-R_{n}) -Q_{n}\zeta_{n-1}][( Q_{n}\zeta_{n-1})^{-1}(I-R_{n})-I]\}^{-1}P_{n} \mathbf{1}\\
 & = & \mathbf{e_{i}}  \{[(I-R_{n}) -Q_{n}\zeta_{n-1}][( Q_{n}\zeta_{n-1})^{-1}](I-R_{n}-Q_{n}\zeta_{n-1})\}^{-1}P_{n} \mathbf{1}\\
 & = & \mathbf{e_{i}}  [(I-Q_{n}\zeta_{n-1}-R_{n})( Q_{n}\zeta_{n-1})^{-1}(I-Q_{n}\zeta_{n-1}-R_{n})]^{-1}P_{n} \mathbf{1}\\
 & = & \mathbf{e_{i}} (I-Q_{n}\zeta_{n-1}-R_{n})^{-1}Q_{n}\zeta_{n-1}(I-Q_{n}\zeta_{n-1}-R_{n})^{-1} P_{n} \mathbf{1}\\
 & = & \mathbf{e_{i}}  (I-Q_{n}\zeta_{n-1}-R_{n})^{-1}Q_{n}\zeta_{n-1}\zeta_{n} \mathbf{1}\\
 & = & \mathbf{e_{i}}  A_{n} \zeta_{n-1}\zeta_{n} \mathbf{1}\\
 & = & \mathbf{e_{i}}  A_{n} \mathbf{1} .
 \eeqnn

\noindent{\large{\bf Acknowledgements:}} The authors would like to
thank  Hongyan Sun, Huaming Wang, Lin
Zhang and Zhou Ke for the stimulating discussions.


{\center

\end{document}